\documentclass[11pt]{article}

\usepackage{amsmath,amsmath,amsthm,amscd,amssymb,verbatim,epsf}
\numberwithin{equation}{section}

\usepackage{amsthm, amssymb, amscd}

\newtheorem{pro}{Proposition}[section]
\newtheorem{thm}[pro]{Theorem}
\newtheorem{lem}[pro]{Lemma}

\def\G{{\Gamma}}
 \def\d{{\delta}}

 \def\L{{\Lambda}}

 \def\a{{\alpha}}
 \def\b{{\beta}}
 \def\p{{\partial}}

 \def\g{{\gamma}}

 \def\2{{\mathbb Z_2}}

 \def\sl2{{SL(2,\mathbb C)}}

 \def\sl{{{\mbox{\tiny $\L$}}}}

\def\wt{\widetilde}
\def\mc{\mathcal}

\def\d{\delta}

\def\a{\alpha}
\def\b{\beta}

\def\L{\Lambda}

\def\g{\gamma}

\begin{document}

\title{Kleinian groups with ubiquitous surface subgroups}
\author{Joseph D. Masters}

\maketitle

\begin{abstract}
We show that every finitely-generated free subgroup of a right-angled,
 co-compact Kleinian reflection group is contained in a surface subgroup.
\end{abstract}

\section{Introduction}

It is conjectured that every co-compact Kleinian group contains
 a surface subgroup. We
 show that, for some special examples, much more is true.

\begin{thm} \label{main}
Let $P$ be a right-angled, compact Coxeter polyhedron in $\mathbb{H}^3$,
 and let $\G(P) \subset Isom(\mathbb{H}^3)$ be the group
 generated by reflections in the faces
 of $P$. Then every finitely-generated
 free subgroup of $\G(P)$ is contained in a surface subgroup
 of $\G(P)$.
\end{thm}

\noindent
\textit{Remarks:}\\
1. It is well-known that every such $\G(P)$ contains a surface
 subgroup.  Indeed, it was shown in \cite{M}
 that the number of ``inequivalent'' surface subgroups
 of $\G(P)$ grows factorially with the genus.\\
\\
2. Lewis Bowen has recently applied Theorem \ref{main}
 to show that every such $\G(P)$ contains a sequence of surface subgroups
 for which the Hausdorff dimensions of the limit sets approach two
 (see \cite{B}).

\section{Outline of the proof}

 Given a free subgroup $F$, we look at the convex
 core $Core(F) = Hull(\Lambda(F))/F$,
 which will be homeomorphic to a handlebody.
 Replacing $Hull(\Lambda(F))$
 with a suitable nighborhood in $\mathbb{H}^3$, we can expand the handlebody
 to make it polyhedral,
 so that the boundary is a union of copies of the
 faces of $F$.  By expanding further, we can make
 the induced decomposition of the boundary finer and finer.
 If we expand enough, it becomes possible
 to attach mirrors to certain faces along
 the boundary (see Figure 1), in such a way that the resulting 3-orbifold
 is the product of a compact 2-orbifold with an interval.
 The desired surface group is a finite-index
 subgroup of the 2-orbifold group.

\begin{figure}[!ht] \label{case3}
{\epsfxsize=3in \centerline{\epsfbox{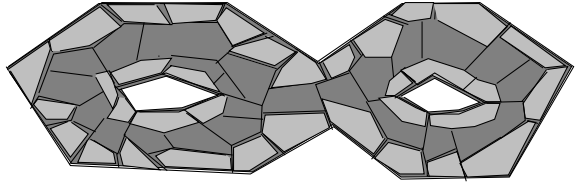}}\hspace{10mm}}
\caption{Mirrors are attached to the lightly-shaded faces.}
\end{figure}

\section{Proof}

\begin{proof}
 The first ingredient is the Tameness Theorem.
 Let $F$ be a free subgroup of $\G(P)$.  Then by \cite{A} and \cite{CG},
 the (infinite-volume) hyperbolic manifold $\mathbb{H}^3/F$ is 
 \textit{topologically tame}, i.e. homeomorphic to the interior
 of a compact 3-manifold.
 Then work of Canary (\cite{C}) implies that
 $F$ is geometrically finite-- i.e. if $C$ is the convex hull
 of the limit set of $F$, then $C/F$ is compact.

 The next step, based on the ideas of \cite{Sc},
 is to give a polyhedral structure to $C/F$.
 Let $\mathcal{T}$ be the tesselation of $\mathbb{H}^3$ by copies of $P$,
 and let $C^+$ be the \textit{tiling hull}
 of $F$-- this is the union of all the tiles in $\mathcal{T}$ which
 meet $C$.
 Then $C^+$ is invariant under $F$, and $C^+/F$ is a compact,
 irreducible 3-manifold with free fundamental group.
 Thus $C^+/F$ is a handlebody $W$.

 The tesselation $\mathcal{T}$ induces
 a tesselation of $\p W$. Since all dihedral angles of
 $P$ are $\pi/2$, then every pair of adjacent faces
 in $\p W$ will meet at an angle of either $\pi/2$ or $\pi$.
 However, if two faces meet at an angle of $\pi$, then we actually
 consider them as part of a single face. 
 Thus, every face in $\p W$ can be decomposed as a union
 $F = X_1 \cup ... \cup X_m$, where each $X_i$ is congruent
 to a face of the original polyhedron $P$.
  Along each $X_i$, we may attach to $W$ a copy of $P$, to obtain a
 handlebody with convex boundary containing $W$, called
 the \textbf{expansion of W along F}.
 More generally, we define an \textbf{expansion} of $W$
 to be a handlebody $W^{\prime} \supset W$, obtained from $W$ by a finite
 sequence of such operations.

 Let $g$ be the genus of $H$, and
 represent $H$ as $P \times I$, for a planar surface $P$.
 Let $\a_1, ..., \a_{g+1}$ be the boundary curves of $P \times \{ 0 \}$.
  Say that a collection of faces $\mc{F}$ of $\p H$ forms a
 \textbf{face annulus}
 if the faces can be indexed $F_1, ..., F_n$, where $F_i$ is adjacent
 to $F_j$ if and only if $|i-j| = 1$ (mod n), and $\cap_i F_i = \emptyset$.
 The last condition excludes the case of three faces meeting at a vertex.

 The following lemma is the key to proving Theorem \ref{main}.

\begin{lem} \label{annuli}
 There is an expansion $W^{\prime}$ of $W$,
 and a collection $\mc{F}$
 of disjoint face annuli $A_1, ..., A_n \subset \p W^{\prime}$,
 so that the core curve of $A_i$ is freely homotopic to $\a_i$ in $W^{\prime}$.
\end{lem}

\begin{proof}
 Let $\mathcal{A} = \cup_i \a_i$.
 Our first claim is that there is an expansion $W^{\prime}$ of $W$
 so that, after an isotopy of the $\a_i$'s to $\p W^{\prime}$,
 we have  $F \cap \mathcal{A}$ being connected for each $F \in \p W^{\prime}$.

 We may assume, after an isotopy,
 that each face in $\p W$ meets $\mathcal{A}$ in a collection
 of disjoint, properly embedded arcs. Let
 $$ k = k(\mathcal{A}) = Max_{F \in \p W} | F \cap \mathcal{A} |.$$
 Suppose $k > 1$.
 Let $n(\mathcal{A})$ be the number of faces in $\p W$
 which meet $\mathcal{A}$ in $k$ components.
 Let $F \in \p W$ such that $|F \cap \mathcal{A}| = k$, and
 let $W^{\prime}$ be the expansion of $W$ along $F$. 
  Note that $W^{\prime} - W$ is a polyhedron
 $P^{\prime}$ (made up of copies of $P$) with dihedral angles $\pi/2$.
 Let $F^{\prime}$ be the face of $P^{\prime}$ which is identified to $F$,
 and let $F_1^{\prime}, ..., F_n^{\prime}$
 be the faces in $P^{\prime}$ which are adjacent
 to $F^{\prime}$, in cyclic order. 

 Let $N_1(F^{\prime})
 = F^{\prime} \cup F_1^{\prime} \cup ... \cup F_n^{\prime}$,
 and let $N_2(F^{\prime})$ be the union of $N_1(F^{\prime})$ together with 
 all faces in $P^{\prime}$ which meet faces in $N_1(F^{\prime})$.
 Since $P^{\prime}$ is a Coxeter polyhedron in $\mathbb{H}^3$,
 it follows that $int \, N_2(F^{\prime})$ is an embedded disk.

 Recall that $\mathcal{A} \cap F$ consists of $k$ disjoint arcs;
 let $\b_1, ..., \b_k$ be the images of these arcs in
 $F^{\prime}$, and
 let $(p_i, q_i)$ be the endpoints of $\b_i$.

\begin{lem} \label{gamma}
 There are disjoint arcs $\g_i$ in $\p P^{\prime} - F^{\prime}$,
 with endpoints $(p_i, q_i)$, so that:\\
1. $|F^* \cap (\cup \g_i)| < k$, for all faces $F^*$ in
 $\p P^{\prime} - N_1(F^{\prime})$.\\
2. $|F_j^{\prime} \cap (\cup \g_i)| = |F_j^{\prime} \cap (\cup \p \b_i)|$
 for all $j$.
\end{lem}

\begin{proof} (Of Lemma \ref{gamma})\\
\\
\textbf{Case 1:}  There are four endpoints (say
 $(p_1, q_1), (p_2, q_2)$) which lie on four distinct sides of $F^{\prime}$.\\

 In this case, we
 let $\d$ be a properly embedded arc in $N_1(F^{\prime})$, disjoint from
 $\cup_i \b_i$,
 which  separates $\b_1$ and $\b_2$ (See Figure 2). 
 For each $i$, let $\b_i^+$ (resp. $\b_i^-$)
 be an arc, properly embedded in some $F_j$, so that one endpoint
 is on $\p N_1(F)$,  the other is the point $p_i$ (resp. $q_i$),
 and so that the arcs $\b_1^{\pm}, \b_2^{\pm}, ... $ are all disjoint
 from each other and from $\d$.
 Let $\b_i^*$ be the component of $\p N_1(F^{\prime}) - (\b_i^+ \cup \b_i^-)$
 which is disjoint from $\d$. 
 Let $\g_i = \b_i^+ \cup \b_i^-  \cup \b_i^*$.
 After an isotopy (supported in a neighborhood of
 $\b_i^*$ in $N_2(F^{\prime}) - int \, N_1(F^{\prime})$) the 
 arcs $\g_i$ satisfy the hypotheses of the lemma.\\
\\
\begin{figure}[!ht] \label{gammafig}
{\epsfxsize=3in \centerline{\epsfbox{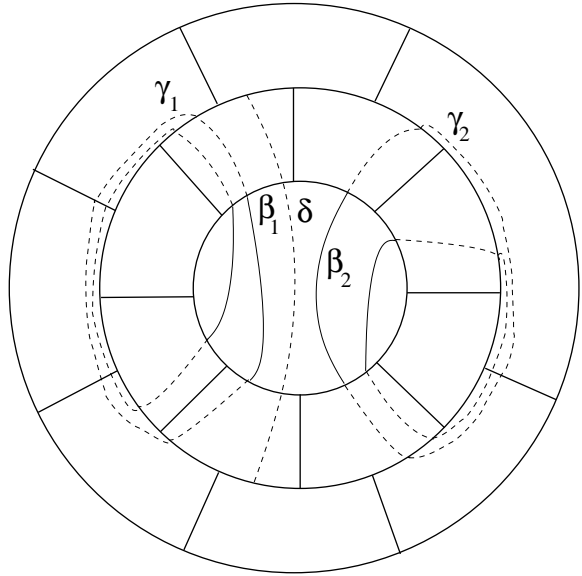}}\hspace{10mm}}
\caption{Construction of $\g_i$'s (Case 1).}
\end{figure}
\textbf{Case 2:}
 Suppose that some edge of $F$ meets every arc $\b_i$.\\

 We repeat the construction
 from Case 1. (i.e.  pick an arc $\d$ in $N_1(F^{\prime})$ disjoint
 from the $\b_i$'s,  separating
 $\b_1$ and $\b_2$; then construct $\b_i^{\pm}$'s, $\b_i^*$'s,
 and $\g_i$'s.)
 The only difference is that we must arrange that the
 arcs $\b_1^+, \b_2^+, ...$ are not all parallel (i.e. their union
 meets at least three distinct sides),
 and that the arcs $\b_1^-, \b_2^-, ...$ are not all parallel.
 This can be done, since, $P^{\prime}$ being a right-angled Coxeter polyhedron
 in $\mathbb{H}^3$, each $F_i^{\prime}$ has at least five edges.
 (See Figure 3).

\begin{figure}[!ht] \label{case2}
{\epsfxsize=3in \centerline{\epsfbox{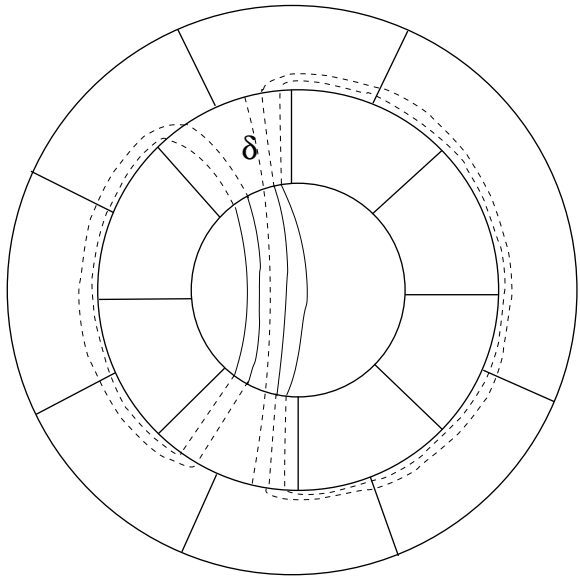}}\hspace{10mm}}
\caption{Construction of $\g_i$'s (Case 2).}
\end{figure}
\end{proof}

  Now we return to the proof of Lemma \ref{annuli}.
 We obtain a loop $\a_i^{\prime}$ in $W^{\prime}$,
 by replacing each $\b_j \subset \a_i$ with $\g_j$.
 Let $\mathcal{A}^{\prime} = \cup \a_i^{\prime}$.
  Since
 the face $F$ has been removed,
 and replaced by faces which meet $\mathcal{A}^{\prime}$ in
 fewer than $k$ components, we have
 $n(\mathcal{A}^{\prime}) < n(\mathcal{A})$.

 Similarly, we see that, by enlarging $W$ repeatedly,
 $n(\mathcal{A})$ can be reduced until it reaches 0.
 By further enlargements, we may assume that $k(\mathcal{A}) = 1$.
  So we may assume that $F \cap \mathcal{A}$ is
 connected for each $F$.
 
 Let $A_i$ be the union of the faces which meet $\a_i$. 
  For each face $F$ in $\cup A_i$, let us define the \textbf{overlap}
 of $F$ by the formula:
$$o(F) = (\textrm{Number of faces in } \cup A_i \textrm{ which are adjacent
 to } F) - 2.$$
 Since the core curve of $A_i$ is essential in $W$,
  no point in $\p W$ meets every face in $A_i$.
 Thus, if $o(F) = 0$ for all $F \in \cup A_i$, then
 the $A_i$'s are the disjoint face annuli  we are looking for.

 Let $F$ be a face in $A_i$, let $F_1$ and $F_2$ be the two
 faces in $A_i$ which are consecutive to $F$, and
 let $e_i = F \cap F_i$.  Let $\g_1$ and $\g_2$
 be the components of $\p F - \{e_1 \cup e_2\}$.
 We say that $F$ is \textbf{good} if one of the $\g_i$'s
 is disjoint from the interior of $\cup A_i$.
\\
\\
 \textbf{Case 3:}  Every face in $\cup A_i$ is good.\\

 Let $F$ be a face in some $A_i$, and let 
 $\b = F \cap (\cup \a_i)$. By previous assumption, $\b$ is connected.
 Let $p$ and $q$ be the endpoints of $\b$.
 As before, let $W^{\prime}$ be the enlargement of $W$ along $F$,
 let $P^{\prime} = W^{\prime} - int \, W$, and let $F^{\prime}$
 be the face of $W^{\prime}$ which is identified to $F$.
   Let $F_1^{\prime}, ..., F_n^{\prime}$ be the faces adjacent to $F^{\prime}$
 in $P^{\prime}$,  labeled consecutively,  so that
 $p \in \p F_1^{\prime}$ and $q \in \p F_i^{\prime}$.
 Since $F$ is good, then we may assume that
 none of the faces $F_1^{\prime}, ..., F_i^{\prime}$ is glued
 to a face in $\cup A_i$.
 
 As in the proof of Lemma \ref{gamma}, we replace $\b$ with
 an appropriate arc $\g \subset \p P^{\prime} - F^{\prime}$.
 In this case, we choose arcs $\b^+$ (resp. $\b^-$)
 from $p$ (resp. $q$) to $\p N_1(F^{\prime})$,
 so that $\b^+$ and $\b^-$ each meet only one face of $\p P^{\prime}$.
 We let $\b^*$ be the component
 of $\p N_1(F^{\prime}) - (\b_1^+ \cup \b_1^-)$
 contained in $F_1^{\prime}, ..., F_i^{\prime}$;
 then we perturb $\b^*$ so that it is a properly embedded
 arc in $N_2(F^{\prime}) - N_1(F^{\prime})$. See Figure 4.

\begin{figure}[!ht] \label{case3}
{\epsfxsize=3in \centerline{\epsfbox{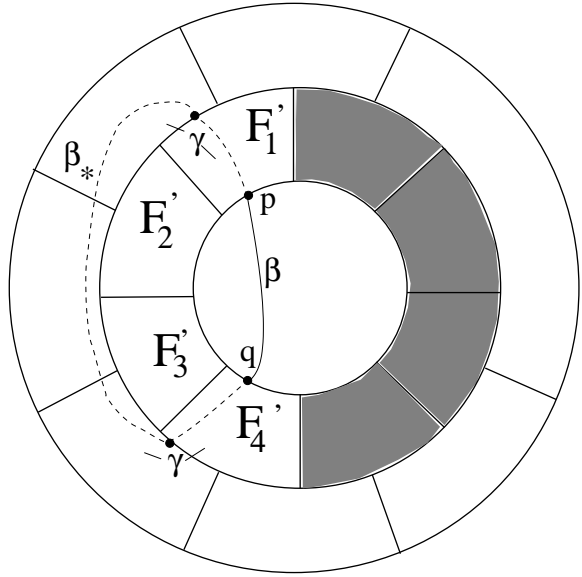}}\hspace{10mm}}
\caption{Construction of $\g$ in Case 3. Only shaded faces can be glued to
 $\cup A_i$.}
\end{figure}

 A complication is that $\p N_2(F^{\prime})$ may not be an embedded
 circle in $P^{\prime}$, and thus there may be pairs of adjacent
 faces in $P^{\prime}$ which meet $\b^*$ non-consecutively.
 In this case, we perform ``shortcut'' operations
 on $\b^*$, as indicated in Figure 5.

\begin{figure}[!ht] \label{case3b}
{\epsfxsize=3in \centerline{\epsfbox{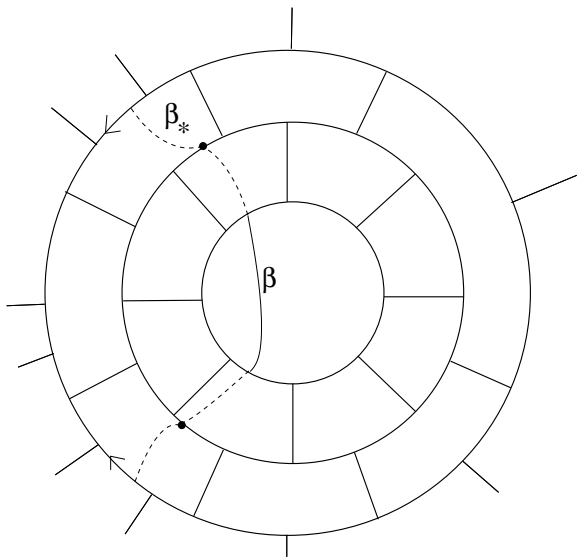}}\hspace{10mm}}
\caption{If the edges with arrows are actually the same,
 then it is possible to shorten the arc $\b^*$.}
\end{figure}

 Let $\g = \b^+ \cup \b^- \cup \b^*$. Then
 we have the required arc $\g$, and a new loop $\a^{\prime}$.
 The number of faces with positive overlap decreases, so 
 eventually we may eliminate them all.\\
 \\
\textbf{Case 4:}  Suppose there is a face $F$ in $\cup A_i$
 which is not good.\\

Here the construction is similar to the construction of Case 3.
 In this case, we choose  $\b^*$  to be 
 either of the two components of $\p N_1(F^{\prime}) - \b_1^+ \cup \b_1^-$;
 then we push $\b^*$ off of $\p N_1(F^{\prime})$; and then,
 as in Case 3, we perform shortcuts if possible.
 The result is that the face $F$ is removed, and replaced
 with good faces.  Repeating this operation along all faces which are
 not good, we may reduce to Case 3.

Thus, we have shown that, after a sequence of enlargements,
 every face in $\cup A_i$ has zero overlap.
 Thus we have constructed the required $A_i$'s, completing
 the proof of Lemma \ref{annuli}.

 \end{proof}

 Returning to the proof of Theorem \ref{main},
 we let $G$ be the group generated by $F$, together with the
 reflections in the lifts to $\mathbb{H}^3$
 of the faces of the face annuli $A_1, ..., A_n$.
 Then we claim that $G$ is the group of a closed, hyperbolic  2-orbifold.

 Indeed, let $V$ be the orbifold with underlying space
 $W$, and with mirrors on the faces of $A_1, ..., A_n$.
 Then $V$ is a hyperbolic 3-orbifold with convex boundary,
 and there is a local isometry $i: V \rightarrow \mathbb{H}^3/\G(P)$,
 with induced map $i_*: \pi_1^{orb}(V) \rightarrow \G(P)$,
 so that $Image(i_*) = G$.
 Since $V$ has convex boundary, every element in
 $\pi_1^{orb}(V)$ is represented by a closed geodesic, and
 since $i$ takes geodesics to geodesics,
 it follows that $i$ is $\pi_1$-injective.

 Note that $V$ is equivalent to a product orbifold
 $X \times I$, where $X$ is the 2-orbifold with reflector
 edges corresponding to one of the components
 of $\p W - \bigcup A_i$.  Thus $G = image(i_*)$ is isomorphic to the orbifold
 fundamental group of $X$.

 The orientable double cover of $X$ is a 2-orbifold, $\wt{X}$,
 where the underlying space is an orientable surface of genus $g$,
 and the cone points of $\wt X$ all have order 2.
 If we identify $G$ with
 $\pi_1^{orb} X$, then the loops generating $F$ all lift to $\wt X$,
 and so $F \subset \pi_1 \wt X$. The group $\pi_1 \wt X$
 has a torsion-free subgroup, of index two (if the number of cone points
 is even) or four (if the number of cone points is odd), containing $F$.
 This is the surface subgroup we were looking for.
  \end{proof}

\noindent
SUNY Buffalo\\
jdmaster@buffalo.edu


\begin{thebibliography}{2}
\bibitem{A} Ian Agol, ``Tameness of hyperbolic 3-manifolds'',  
 arXiv:math/0405568.\\

\bibitem{B} L. Bowen, ``Free groups in lattices'', pre-print.\\

\bibitem{C} R. D. Canary, ``A covering theorem for
 hyperbolic $3$-manifolds and its applications'',
  \textit{Topology}  \textbf{35}  (1996),  no. 3, 751--778.\\

\bibitem{CG}  D. Calegari and D. Gabai,
 ``Shrinkwrapping and the taming of hyperbolic 3-manifolds'',
  \textit{J. Amer. Math. Soc.}  \textbf{19}  (2006),  no. 2, 385--446
 (electronic).

\bibitem{M} J. D. Masters, 
 ``Counting immersed surfaces in hyperbolic 3-manifolds'',
\textit{Algebr. Geom. Topol.} \textbf{5} (2005), 835--864 (electronic). 

\bibitem{Sc}  G. P. Scott, ``Subgroups of surface groups are almost
 geometric'',  \textit{J. London Math. Soc.} (2)  \textbf{17}
  (1978), no. 3, 555--565.
\end{thebibliography}
\end{document}